\documentclass[12pt]{article}
\usepackage{latexsym}
\usepackage{amsmath}
\usepackage{amssymb}
\usepackage{amsfonts}
\usepackage[active]{srcltx}

\numberwithin{equation}{section}

\newtheorem{theorem}{Theorem}[section]
\newtheorem{lemma}[theorem]{Lemma}

\newtheorem{remark}[theorem]{Remark}
\newcommand{\eproof}{{\mbox{\ }~\hfill
\mbox{\large $\Box$} \par \vskip 10pt}}

\newcommand{\R}{{\mathbb R}}

\renewcommand{\div}{{\rm div}}

\title{Quantitative strong unique continuation for the Lam\'e system with less regular coefficients}
\author{C-L Lin\thanks{Department of Mathematics, NCTS, National Cheng Kung University,
Tainan 701, Taiwan. Email:cllin2@mail.ncku.edu.tw}  \quad G
Nakamura\thanks{Department of Mathematics, Hokkaido University,
Sapporo 060-0810, Japan.(Email: gnaka@math.sci.hokudai.ac.jp)}\quad G
Uhlmann\thanks{Department of Mathematics, University of Washington,
Box 354350, Seattle 98195-4350, USA.
Email:gunther@math.washington.edu}\quad J-N Wang\thanks{Department
of Mathematics, Taida Institute of Mathematical Sciences, NCTS
(Taipei), National Taiwan University, Taipei 106, Taiwan.
Email:jnwang@math.ntu.edu.tw}}

\date{}

\begin{document}
\maketitle

\begin{abstract}
In this paper we prove a quantitative form of the strong unique
continuation property for the Lam\'e system when the Lam\'e
coefficients $\mu$ is Lipschitz and
$\lambda$ is essentially bounded in dimension $n\ge 2$. This result
is an improvement of our earlier result \cite{lin5} in which both
$\mu$ and $\lambda$ were assumed to be Lipschitz.
\end{abstract}

\section{Introduction}\label{sec1}
\setcounter{equation}{0}

Assume that $\Omega$ is a connected open set containing $0$ in
$\R^n$ for $n\geq 2$. Let $\mu(x)\in C^{0,1}(\Omega)$ and
$\lambda(x),\rho(x)\in L^{\infty}(\Omega)$ satisfy
\begin{equation}\label{1.1}
\begin{cases}
\mu(x)\geq\delta_0,\quad\quad
\lambda(x)+2\mu(x)\geq\delta_0\quad\forall\
\text{a.e.}\ x\in\Omega,\\
\|\mu\|_{C^{0,1}(\Omega)}+\|\lambda\|_{L^{\infty}(\Omega)}\leq
M_0,\quad \|\rho\|_{L^{\infty}(\Omega)}\le M_0
\end{cases}
\end{equation}
with positive constants $\delta_0, M_0$, where we define
$$\|f\|_{C^{0,1}(\Omega)}=\|f\|_{L^{\infty}(\Omega)}+\|\nabla
f\|_{L^{\infty}(\Omega)}.$$ The isotropic elasticity system, which
represents the displacement equation of equilibrium, is given by
\begin{equation}\label{1.2}
\text{div}(\mu(\nabla u+(\nabla
u)^t))+\nabla(\lambda\text{div}u)+\rho u=0\quad\text{in}\ \Omega,
\end{equation}
where $u=(u_1,u_2,\cdots,u_n)^t$ is the displacement vector and $(\nabla
u)_{jk}=\partial_ku_j$ for $j,k=1,2,\cdots,n$.

We are interested in the strong unique continuation property (SUCP)
of \eqref{1.2}. More precisely, we would like to show that any
nontrivial solution of \eqref{1.2} can only vanish of finite order at any point of
$\Omega$. We also give an estimate of the vanishing
order for $u$, which can be seen as a quantitative description of
the SUCP for \eqref{1.2}. Here we list some of the known results on the SUCP for
\eqref{1.2}:

\begin{itemize}
\item $\lambda,\mu\in C^{1,1}$, $n\ge 2$ (quantitative): Alessandrini and Morassi
\cite{almo}.

\item $\lambda,\mu\in C^{0,1}$, $n=2$ (qualitative): Lin and Wang \cite{lw05}.

\item $\lambda\in L^{\infty},\mu\in C^{0,1}$, $n=2$ (qualitative):
Escauriaza \cite{es}.

\item $\lambda,\mu\in C^{0,1}$, $n\ge 2$ (quantitative): Lin,
Nakamura, and Wang \cite{lin5}.
\end{itemize}

In this paper, we relax the regularity assumption on $\lambda$ in
\cite{lin5} to $\lambda\in L^{\infty}(\Omega)$. In view of
counterexamples by Plis \cite{pl} or Miller \cite{mi}, this regularity assumption seems
to be optimal. This improvement was inspired by our recent work on
the Stokes system \cite{lin6}. We now state the main results of the
paper. Assume that there exists $0<R_0\le 1$ such that
$B_{R_0}\subset\Omega$. Hereafter $B_r$ denotes an open ball of
radius $r>0$ centered at the origin.
\begin{theorem}\rm{(Optimal three-ball inequalities)}\label{thm1.1}
There exists a positive number $\tilde{R}<1$, depending only on
$n,M_0,\delta_0$, such that if $\ 0<R_1<R_2<R_3\leq R_0$ and
$R_1/R_3<R_2/R_3<\tilde{R}$, then
\begin{equation}\label{1.7}
\int_{|x|<R_2}|u|^2dx\leq
{C}\left(\int_{|x|<R_1}|u|^2dx\right)^{\tau}\left(\int_{|x|<{R_3}}|u|^2dx\right)^{1-\tau}
\end{equation}
for $u\in H_{loc}^1({B}_{R_0})$ satisfying \eqref{1.2} in
${B}_{R_0}$, where the constant ${C}$ depends on $R_2/R_3$, $n$,
$M_0,\delta_0$, and $0<\tau<1$ depends on $R_1/R_3$, $R_2/R_3$,
$n,M_0,\delta_0$. Moreover, for fixed $R_2$ and $R_3$, the exponent
$\tau$ behaves like $1/(-\log R_1)$ when $R_1$ is sufficiently
small.
\end{theorem}

\begin{theorem}\label{thm1.2}
Let $u\in H^1_{loc}(\Omega)$ be a nontrivial solution of \eqref{1.2}, then
there exist positive constants $K$ and $m$, depending on
$n,M_0,\delta_0$ and $u$, such that
\begin{equation}\label{1.8}
\int_{|x|<R}|u|^2 dx\ge KR^m
\end{equation}
for all $R$ sufficiently small.
\end{theorem}
\begin{remark}\label{rem1.2}
Based on Theorem~\ref{thm1.1}, the constants $K$ and $m$ in
\eqref{1.4} are explicitly given by
$$
K=\int_{|x|<R_3}|u|^2dx
$$
and
$$
m=\tilde
C+\log\Big{(}\frac{\int_{|x|<R_3}|u|^2dx}{\int_{|x|<R_2}|u|^2dx}\Big{)},
$$
where $\tilde C$ is a positive constant depending on
$n,M_0,\delta_0$ and $R_2/R_3$.
\end{remark}

\section{Reduced system and estimates}\label{sec2}
\setcounter{equation}{0}

Here we want to find a reduced system from \eqref{1.2}. This is a crucial step in our approach. Let us write \eqref{1.2} into a non-divergence form:
\begin{equation}\label{1.3}
\mu\Delta u+\nabla((\lambda+\mu)\ \div u)+(\nabla u+(\nabla u)^t)\nabla\mu-\div u\nabla\mu+\rho u =0.
\end{equation}
Dividing \eqref{1.3} by $\mu$ yields
\begin{eqnarray}\label{1.4}
&&\Delta u+\frac{1}{\mu}\nabla((\lambda+\mu)\ \div u)+(\nabla u+(\nabla u)^t)\frac{\nabla\mu}{\mu}-\div u\frac{\nabla\mu}{\mu}+\frac{\rho}{\mu} u\notag\\
&=&\Delta u+\nabla(\frac{\lambda+\mu}{\mu}\ \div u)+(\nabla u+(\nabla u)^t)\frac{\nabla\mu}{\mu}-\div u(\frac{\nabla\mu}{\mu}+(\lambda+\mu)\nabla(\frac{1}{\mu}))\notag\\
&&+\frac{\rho}{\mu} u\notag\\
&=&\Delta u+\nabla(a(x)v)+G\notag\\
&=&0,
\end{eqnarray}
where $$a(x)=\frac{\lambda+\mu}{\lambda+2\mu}\in L^{\infty}(\Omega),\quad v=\frac{\lambda+2\mu}{\mu}\ \div u$$ and
$$G=(\nabla u+(\nabla u)^t)\frac{\nabla\mu}{\mu}-\div u(\frac{\nabla\mu}{\mu}+(\lambda+\mu)\nabla(\frac{1}{\mu}))+\frac{\rho}{\mu} u.$$
Taking the divergence on \eqref{1.4} gives
\begin{equation}\label{1.5}
\Delta v+\div G=0.
\end{equation}
Our reduced system now consists of \eqref{1.4} and \eqref{1.5}. It follows easily from \eqref{1.5} that if $u\in H^1_{loc}(\Omega)$, then $v\in H^1_{loc}(\Omega)$.

To prove the main results, we rely on suitable Carleman estimates.
Denote $\varphi_{\beta}=\varphi_{\beta}(x) =\exp
(-\beta\tilde{\psi}(x))$, where $\beta>0$ and $\tilde{\psi}(x)=\log
|x|+\log((\log |x|)^2)$. Note that $\varphi_{\beta}$ is less
singular than $|x|^{-\beta}$. We use the notation $X\lesssim Y$ or
$X\gtrsim Y$ to mean that $X\le CY$ or $X\ge CY$ with some constant
$C$ depending only on $n$.

\begin{lemma}{\rm\cite[Lemma 2.4]{lin5}}\label{lem2.1}
There exist a sufficiently small number $r_1>0$ depending on $n$ and
a sufficiently large number $\beta_1>3$ depending on $n$ such that
for all $w\in U_{r_1}$ and $f=(f_1,\cdots,f_n)\in (U_{r_1})^{n}$,
$\beta\geq \beta_1$, we have that
\begin{eqnarray}\label{2.1}
&&\int\varphi^2_\beta (\log|x|)^2(\beta|x|^{4-n}|\nabla w|^2+\beta^3|x|^{2-n}|w|^2)dx\notag\\
&\lesssim& \int \varphi^2_\beta
(\log|x|)^{4}|x|^{2-n}[(|x|^{2}\Delta w+|x|{\rm div}
f)^2+\beta^2\|f\|^2]dx,
\end{eqnarray}
where $U_{r_1}=\{w\in C_0^{\infty}(\R^n\setminus\{0\}): \mbox{\rm
supp}(w)\subset B_{r_0}\}$.
\end{lemma}
Next, replacing $\beta$ by $\beta+1$ in \eqref{2.1}, we get another Carleman estimate.
\begin{lemma}\label{lem2.2}
There exist a sufficiently small number $r_1>0$ depending on $n$ and
a sufficiently large number $\beta_1>2$ depending on $n$ such that
for all $w\in U_{r_1}$ and $f=(f_1,\cdots,f_n)\in (U_{r_1})^{n}$,
$\beta\geq \beta_1$, we have that
\begin{eqnarray}\label{2.2}
&&\int\varphi^2_\beta (\log|x|)^{-2}(\beta|x|^{2-n}|\nabla w|^2+\beta^3|x|^{-n}|w|^2)dx\notag\\
&\lesssim& \int \varphi^2_\beta
|x|^{-n}[(|x|^{2}\Delta w+|x|{\rm div}
f)^2+\beta^2\|f\|^2]dx.
\end{eqnarray}
\end{lemma}

In addition to Carleman estimates, we also need the following
Caccioppoli's type inequality.
\begin{lemma}\label{lem3.1}
Let $u\in (H^1_{loc}(\Omega))^{n}$ be a solution of
\eqref{1.1}. Then for any $0<a_3<a_1<a_2<a_4$ such that
$B_{a_4r}\subset\Omega$ and $|a_4r|<1$, we have
\begin{equation}\label{3.1}
\int_{a_1r<|x|<a_2r}|x|^{4}|\nabla v|^2+|x|^{2}|v|^2+|x|^{2}|\nabla
u|^2dx\le C_0\int_{a_3r<|x|<a_4r}|u|^2dx
\end{equation}
where the constant $C_0$ is independent of $r$ and $u$. Here
$v$ is defined in \eqref{1.4}.
\end{lemma}
The proof of Lemma~\ref{lem3.1} will be given in the next section. Here we would like to outline how to proceed the proofs of main theorems. The detailed arguments can be found in \cite{lin5} or \cite{lin6}. Firstly, applying \eqref{2.2} to $w=u$, $f=|x|a(x)v$ and using \eqref{1.4}, we have that
\begin{eqnarray}\label{2.3}
&&\int\varphi^2_\beta (\log|x|)^{-2}(\beta|x|^{2-n}|\nabla u|^2+\beta^3|x|^{-n}|u|^2)dx\notag\\
&\lesssim& \int \varphi^2_\beta|x|^{-n}[\big{(}|x|^{2}\Delta u+|x|{\rm div}(|x|a(x)v)\big{)}^2+\beta^2\||x|a(x)v\|^2]dx.
\end{eqnarray}
Next, applying \eqref{2.1} to $w=v$, $f=|x|G$ and using \eqref{1.5}, we get that
\begin{eqnarray}\label{2.4}
&&\int\varphi^2_\beta (\log|x|)^2(\beta|x|^{4-n}|\nabla v|^2+\beta^3|x|^{2-n}|v|^2)dx\notag\\
&\lesssim& \int \varphi^2_\beta
(\log|x|)^{4}|x|^{2-n}[\big{(}|x|^{2}\Delta v+|x|{\rm div}(|x|G)\big{)}^2+\beta^2\||x|G\|^2]dx.\notag\\
\end{eqnarray}
Finally, adding $\beta\times$\eqref{2.3} and \eqref{2.4} together and using \eqref{3.1}, we can prove Theorem~\ref{thm1.1} and \ref{thm1.2} as in \cite{lin5} and \cite{lin6}.

\section{Proof of Lemma~\ref{lem3.1}}\label{sec3}
\setcounter{equation}{0}

Define $b_1=(a_1+a_3)/2$ and $b_2=(a_2+a_4)/2$. Let
$X=B_{a_4r}\backslash \bar{B}_{a_3r}$, $Y=B_{b_2r}\backslash
\bar{B}_{b_1r}$ and $Z=B_{a_2r}\backslash \bar{B}_{a_1r}$. Let
$\xi(x)\in C^{\infty}_0 ({\mathbb R}^n)$ satisfy $0\le\xi(x)\leq 1$
and
\begin{eqnarray}\label{3.2}
\xi (x)=
\begin{cases}
\begin{array}{l}
0,\quad |x|\leq a_3r,\\
1,\quad b_1r<|x|<b_2r,\\
0,\quad |x|\geq a_4r.
\end{array}
\end{cases}
\end{eqnarray}
From \eqref{1.2}, we have that
\begin{eqnarray}\label{3.3}
0&=&-\int[\text{div}(\mu(\nabla u+(\nabla u)^t))+\nabla(\lambda\text{div}u)+\rho u]\cdot (\xi^2\bar{u})dx\notag\\
&=&\int\sum_{ijkl=1}^n[\lambda\delta_{ij}\delta_{kl}+\mu(\delta_{il}\delta_{jk}+\delta_{ik}\delta_{jl})]\partial_{x_l}u_k\partial_{x_j}(\xi^2\bar{u}_i)dx-\int\rho\xi^2|u|^2dx\notag\\
&=&\int\xi^2\sum_{ijkl=1}^n[\lambda\delta_{ij}\delta_{kl}+\mu(\delta_{il}\delta_{jk}+\delta_{ik}\delta_{jl})]\partial_{x_l}u_k\partial_{x_j}\bar{u}_idx\notag\\
&&+\int\sum_{ijkl=1}^n\partial_{x_j}(\xi^2)[\lambda\delta_{ij}\delta_{kl}+\mu(\delta_{il}\delta_{jk}+\delta_{ik}\delta_{jl})]\partial_{x_l}u_k\bar{u}_idx-\int\rho\xi^2|u|^2dx\notag\\
&=&I_1+I_2,
\end{eqnarray}
where
$$I_1=\int\xi^2[\sum_{ij=1}^n\lambda\partial_{x_j}u_j\partial_{x_i}\bar{u}_i+\sum_{ij=1}^n\mu(\partial_{x_i}u_j\partial_{x_j}\bar{u}_i+\partial_{x_j}u_i\partial_{x_j}\bar{u}_i)]dx$$
and
$$I_2=\int\sum_{ijkl=1}^n\partial_{x_j}(\xi^2)[\lambda\delta_{ij}\delta_{kl}+\mu(\delta_{il}\delta_{jk}+\delta_{ik}\delta_{jl})]\partial_{x_l}u_k\bar{u}_idx-\int\rho\xi^2|u|^2dx.$$

Observe that
\begin{eqnarray}\label{3.4}
&&\int\xi^2(2\mu-\frac{\delta_0}{2}) \partial_{x_i}u_j\partial_{x_j}\bar{u}_idx\notag\\
&=& -\int\partial_{x_j}[\xi^2(2\mu-\frac{\delta_0}{2})]\partial_{x_i}u_j\bar{u}_idx-\int\xi^2(2\mu-\frac{\delta_0}{2}) \partial^2_{x_ix_j}u_j\bar{u}_idx\notag\\
&=& -\int\partial_{x_j}[\xi^2(2\mu-\frac{\delta_0}{2})]\partial_{x_i}u_j\bar{u}_idx+\int\partial_{x_i}[\xi^2(2\mu-\frac{\delta_0}{2})]\partial_{x_j}u_j\bar{u}_idx\notag\\
&&+\int\xi^2(2\mu-\frac{\delta_0}{2})
\partial_{x_j}u_j\partial_{x_i}\bar{u}_idx.
\end{eqnarray}
It follows from \eqref{3.4} that
\begin{eqnarray}\label{3.5}
I_1&=&\int\xi^2[\sum_{ij=1}^n\lambda\partial_{x_j}u_j\partial_{x_i}\bar{u}_i+\sum_{ij=1}^n(2\mu-\frac{\delta_0}{2})(\partial_{x_i}u_j\partial_{x_j}\bar{u}_i)]dx\notag\\
&&+\int\sum_{ij=1}^n\xi^2(\mu-\frac{\delta_0}{2})(\partial_{x_j}u_i\partial_{x_j}\bar{u}_i-\partial_{x_i}u_j\partial_{x_j}\bar{u}_i)dx\notag\\
&&+\frac{\delta_0}{2}\int\sum_{ij=1}^n\xi^2\partial_{x_j}u_i\partial_{x_j}\bar{u}_idx\notag\\
&=&\int(2\mu+\lambda-\frac{\delta_0}{2})\xi^2\sum_{ij=1}^n(\partial_{x_j}u_j\partial_{x_i}\bar{u}_i)dx\notag\\
&&+\int\sum_{ij=1}^n\xi^2(\mu-\frac{\delta_0}{2})(\partial_{x_j}u_i\partial_{x_j}\bar{u}_i-\partial_{x_i}u_j\partial_{x_j}\bar{u}_i)dx\notag\\
&&+\frac{\delta_0}{2}\int\sum_{ij=1}^n\xi^2\partial_{x_j}u_i\partial_{x_j}\bar{u}_idx+I_3,
\end{eqnarray}
where
$$I_3=\sum_{ij=1}^n\int\partial_{x_i}[\xi^2(2\mu-\frac{\delta_0}{2})]\partial_{x_j}u_j\bar{u}_i-\partial_{x_j}[\xi^2(2\mu-\frac{\delta_0}{2})]\partial_{x_i}u_j\bar{u}_idx.$$
Since
\begin{eqnarray*}
&&\int\sum_{ij=1}^n\xi^2(\mu-\frac{\delta_0}{2})(\partial_{x_j}u_i\partial_{x_j}\bar{u}_i-\partial_{x_i}u_j\partial_{x_j}\bar{u}_i)dx\\
&=&\frac{1}{2}\int\sum_{ij=1}^n\xi^2(\mu-\frac{\delta_0}{2})|\partial_{x_j}u_i-\partial_{x_i}u_j|^2dx,
\end{eqnarray*}
we obtain that
\begin{eqnarray}\label{3.7}
I_1\geq \frac{\delta_0}{2}\int|\xi\nabla u|^2dx+I_3.
\end{eqnarray}
Combining \eqref{3.3} and \eqref{3.7}, we have that
\begin{equation*}
\int_{Y}|\nabla u|^2dx \le\int_{X} |\xi\nabla u|^2 dx\leq
C_1\int_{X}|x|^{-2}|u|^2dx,
\end{equation*}
which implies
\begin{equation}\label{3.8}
\int_{Y}|x|^{2}|\nabla u|^2dx\leq C_2\int_{X}|u|^2dx.
\end{equation}
Here and below all constants $C_1,C_2,\cdots$ depend on $\delta_0$,
$M_0$.

To estimate $\nabla v$, we define $\chi(x)\in C^{\infty}_0 ({\mathbb
R}^n)$ satisfy $0\le\chi(x)\leq 1$ and
\begin{eqnarray*}
\chi (x)=
\begin{cases}
\begin{array}{l}
0,\quad |x|\leq b_1r,\\
1,\quad a_1r<|x|<a_2r,\\
0,\quad |x|\geq b_2r.
\end{array}
\end{cases}
\end{eqnarray*}
By \eqref{1.5}, we derive that
\begin{eqnarray}\label{3.10}
&& \int|\chi(x)\nabla v|^2dx\notag\\
&=&\int\nabla v\cdot\nabla(\chi^2 \bar{v})dx-2\int\chi\nabla v\cdot \bar{v}\nabla\chi dx\notag\\
&\le&|\int(\div G) \chi^2\bar{v} dx|+2\int|\chi\nabla v\cdot \bar{v}\nabla\chi |dx\notag\\
&\le&|\int(\div G) \chi^2\bar{v} dx|+\frac{1}{4}\int|\chi\nabla v|^2dx+C_3\int_{Y}|x|^{-2}|v|^2dx\notag\\
&\le&C_4\int_{Y} |\nabla u|^2 dx+C_4\int_{Y} |u|^2 dx+\frac{1}{2}\int|\chi\nabla v|^2dx+C_4\int_{Y}|x|^{-2}|v|^2dx\notag\\
&\le&C_5\int_{Y} |x|^{-2}|\nabla u|^2 dx+C_4\int_{Y} |u|^2
dx+\frac{1}{2}\int|\chi\nabla v|^2dx.
\end{eqnarray}
Therefore, we get from \eqref{3.10} that
\begin{equation*}
\int_{Z}|\nabla v|^2dx\le2C_5\int_{Y} |x|^{-2}|\nabla u|^2
dx+2C_4\int_{Y} |u|^2 dx
\end{equation*}
and hence
\begin{equation}\label{3.12}
\int_{Z}|x|^{4}|\nabla v|^2dx\le C_6\int_{Y} |x|^{2}|\nabla u|^2
dx+C_6\int_{Y}|x|^{4} |u|^2 dx.
\end{equation}

Putting together $K\times$\eqref{3.8} and \eqref{3.12}, we have that
\begin{eqnarray}\label{3.13}
&&K\int_{Y}|x|^{2}|\nabla u|^2dx+\int_{Z}|x|^{4}|\nabla v|^2dx\notag\\
&\le&KC_2\int_{X}|u|^2dx+C_6\int_{Y} |x|^{2}|\nabla u|^2 dx+C_6\int_{Y}|x|^{4} |u|^2 dx.
\end{eqnarray}
Choosing $K=2C_6$ in \eqref{3.13} yields
\begin{eqnarray*}
&&\int_{Z}|x|^{2}|v|^2dx+\int_{Z}|x|^{2}|\nabla u|^2dx+\int_{Z}|x|^{4}|\nabla v|^2dx\notag\\
&\le&C_7\int_{Y}|x|^{2}|\nabla u|^2dx+C_7\int_{Z}|x|^{4}|\nabla v|^2dx\notag\\
&\le&C_8\int_{X}|u|^2dx,
\end{eqnarray*}
The proof is now complete.\eproof

\section*{Acknowledgements}
Lin was supported in part by the National Science Council of Taiwan.
Nakamura was partially supported by Grant-in-Aid for Scientific
Research (B)(No.19340028) of Japan Society for Promotion of Science.
Uhlmann was supported in part by NSF and a Walker Family Endowed
Professorship. Wang was partially supported by the National Science
Council of Taiwan.

\end{document}